\documentclass[letterpaper, 10 pt, conference]{ieeeconf}
\IEEEoverridecommandlockouts
\usepackage{graphicx}
\usepackage{cite}
\usepackage{amsmath}
\usepackage{amssymb}
\usepackage{mathtools}
\usepackage{algorithm}
\usepackage{algpseudocode}
\usepackage{booktabs}
\usepackage{multirow}
\usepackage[dvipsnames,table]{xcolor}
\usepackage{caption}
\usepackage{subcaption}
\usepackage{siunitx}
\usepackage{tikz}
\usepackage{balance}
\usepackage{pifont}
\usepackage[hidelinks]{hyperref}

\newtheorem{remark}{Remark}

\newcommand{\cupiqp}{\textsc{cuPIQP}}

\newif\ifanonymous
\anonymousfalse
\ifanonymous
\newcommand{\codefootnote}{Code will be released under an open-source license upon acceptance.}
\else
\newcommand{\codefootnote}{\url{https://github.com/PREDICT-EPFL/cupiqp}, BSD-2-Clause license.}
\fi

\title{\LARGE \bf
CuPIQP: A Batched and Differentiable QP Solver on GPUs \\
for Robot Learning and Control
}

\ifanonymous
\author{Anonymous Authors%
\thanks{Submitted for double-anonymous review.}%
}
\else
\author{Fenglong Song$^{1}$, Roland Schwan$^{2}$, Yuwen Chen$^{1}$, Yuning Jiang$^{3}$ and Colin N. Jones$^{1}$%
\thanks{This work was supported by the Swiss National Science Foundation under the NCCR Automation project, grant agreement 51NF40\_180545.}%
\thanks{$^{1}$F. Song, Y. Chen and C. N. Jones are with the Automatic Control Laboratory, EPFL, Switzerland.
{\tt\small \{fenglong.song, yuwen.chen, colin.jones\}@epfl.ch}}%
\thanks{$^{2}$R. Schwan is with the Accessible and Accelerated Robotics Lab, Dartmouth College, USA.
{\tt\small roland.schwan@dartmouth.edu}}%
\thanks{$^{3}$Y. Jiang is with the Institute for Transport Planning and Systems, ETH Z\"urich, Switzerland.
{\tt\small yuning.jiang@ieee.org}}%
}
\fi

\begin{document}

\maketitle

\thispagestyle{empty}
\pagestyle{empty}

\begin{abstract}
Quadratic programs (QPs) are fundamental building blocks of many robotic tasks such as inverse kinematics, whole-body control, and model predictive control. Extending these methods to large-scale robot learning and control requires solving many QPs concurrently and, for end-to-end learning, differentiating through their solutions. Efficient batched optimization is therefore essential to bringing established model-based methods into GPU-based robotics pipelines. We present \cupiqp{}, a batched and differentiable QP solver on GPUs based on a proximal interior-point method, with unified support for dense, sparse, and multistage problems. The solver exploits parallelism across independent QPs and their internal structure while preserving per-problem numerical updates and termination decisions. Benchmarks with up to 16384 QPs per batch demonstrate solve-time speedups of up to $18.8\times$ over \textsc{Moreau} and $4.26\times$ over eight-core \textsc{PIQP}. In an end-to-end learning task, \cupiqp{} accelerates training by $4.6\times$ over \textsc{QPAX} and $7.8\times$ over \textsc{Moreau}. 
By reducing the computational cost, \cupiqp{} provides efficient, reusable infrastructure for integrating model-based optimization into large-scale robot learning and control.
\end{abstract}


\section{Introduction}
\label{sec:intro}

Quadratic programs lie on the critical path of many robot
control tasks.
QP solvers have traditionally been designed to minimize the online latency for solving one or a tiny batch of problems. Existing solvers such as \textsc{qpOASES}~\cite{qpoases}, \textsc{OSQP}~\cite{osqp},
\textsc{HPIPM}~\cite{frison2020hpipm}, \textsc{ProxQP}~\cite{bambade2025proxqp} and \textsc{PIQP}~\cite{piqp} are highly effective under this execution model. 

Nonetheless, an increasingly important class of robot learning and control workloads follows a different execution model, where optimization is embedded in an offline training pipeline and evaluated concurrently over a large number of instances~\cite{amos2018differentiablempc, romero2026actorcritic, reiter2026synthesis, jeon2026residualmpc}. 
In such workloads, the same optimization layer is applied across problem instances arising from parallel environments or sampled trajectories~\cite{rudin2022massively}. The resulting computational primitive is a batch of $B$ independent QPs $\{\mathcal{P}_j\}_{j=1}^B$, where the $j$-th problem
$\mathcal{P}_j$ is given by
\begin{equation}
\label{eq:general_batch_qp}
\begin{aligned}
\min_{x_j}\quad
& \tfrac{1}{2}x_j^\top P_j x_j + c_j^\top x_j, \\
\text{s.t.}\quad
& A_j x_j = b_j, \\
& G_j x_j \leq h_j,
\end{aligned}
\end{equation}
where $P_j\in\mathbb{S}_+^{n}$,
$A_j\in\mathbb{R}^{p\times n}$, and $G_j\in\mathbb{R}^{m\times n}$. 
The dimensions $n$, $p$, and $m$, as well as the sparsity patterns of $P_j$, $A_j$, and $G_j$, are shared across the batch, whereas the numerical values of the problem data may vary across instances. 
Consequently, the relevant objective is not only low latency for an individual QP, but also high throughput over the batch while preserving the solution and termination status of each problem.
This highly concurrent workload is naturally well suited to execution on Graphics Processing Units (GPUs) because the independence of the QPs exposes abundant problem-level parallelism.
The case for a GPU-native solver is further strengthened by data locality. In robot learning pipelines, QP data are often produced on the GPU by a policy, dynamics model, or simulator, and the solutions are immediately consumed there. 
A CPU solver introduces device--host transfers and host-side synchronization, interrupting device-resident
execution and precluding end-to-end CUDA Graph capture across the solver boundary. 
A GPU-native batched solver instead keeps the complete optimization workflow on the device, avoiding these transfers
and synchronization boundaries.


\subsection{Related Work}
\label{subsec:related}

\subsubsection{GPU solvers for single large sparse problems}
A class of GPU-enabled solvers targets individual large-scale sparse optimization problems~\cite{CuClarabel, chari2026gpu, cuopt} where sparse linear algebra is sufficiently intensive to amortize GPU kernel-launch and synchronization overheads, and reducing such overheads is therefore not a primary design objective.
By contrast, these overheads can dominate when solving large batches of small-scale QPs.
Moreover, they do not explicitly support batched solving. 
Although independent QPs can be assembled into a single block-diagonal problem, this formulation represents repeated structure as generic sparsity rather than exploiting it through batch-oriented operations, potentially resulting in suboptimal performance.

\subsubsection{Batched solvers for general QP}
Another class of solvers treats the batch as an explicit computational dimension for general QPs.
\textsc{QPTH}~\cite{amos2017optnet} and \textsc{QPAX}~\cite{tracy2024differentiability} implement a differentiable interior-point method that solves batches of dense QPs concurrently on GPUs.
\textsc{Moreau}~\cite{moreau2026} is a differentiable solver for sparse conic optimization problems with shared sparsity pattern, providing both CPU and GPU backends.
However, each solver is specialized for a single problem structure and does not provide unified support for dense, sparse, and multistage QPs.

\subsubsection{Batched solvers for MPC}
A complementary line of work exploits batching at the level of the complete nonlinear MPC pipeline.
Recent GPU-based frameworks provide batched and differentiable sequential quadratic programming (SQP) implementations for large-batch MPC~\cite{du2026gatogpuacceleratedbatchedtrajectory, adabag2025differentiablemodelpredictivecontrol, bravopalacios2026turbompcfastscalabledifferentiable, hose2026warpmpclargebatchmpcgpu}. These frameworks generate a batch of structurally homogeneous QP subproblems at each SQP iteration, making an efficient batched QP solver a core primitive of the overall pipeline.
On CPUs, \textsc{acados}~\cite{acados,frey2025differentiable} provides a batched and differentiable interface by parallelizing independent solves across multiple processor cores.



\subsection{Contributions}
\label{subsec:contributions}

We present \textsc{cuPIQP} as an open-source software\footnote{\codefootnote}.
Its main technical contributions are summarized as follows:

\begin{itemize}
\item {Batch-native GPU execution.}
We design the solver around batches of independent QPs, exploiting problem-level parallelism for high throughput while retaining independent iterates, convergence decisions, solutions, and termination statuses for all problems in the batch.

\item {Structure-specialized QP backends.}
\textsc{cuPIQP} supports dense, sparse, and multistage QPs through specialized data layouts and linear algebra backends under a unified solver interface, allowing users to exploit the structure of a variety of problem classes.

\item {GPU-resident integration and differentiation.}
Problem data, solver iterates, per-problem statuses, and adjoint computations remain on the GPU. The solver interoperates directly with CUDA-resident arrays without host staging, provides implicit differentiation through a PyTorch interface, and supports CUDA Graph capture on compatible execution paths to reduce CPU orchestration and kernel-launch overhead.
\end{itemize}

\emph{Notation:} $\mathbb{S}^n_+$ and $\mathbb{S}^n_{++}$ denote symmetric positive
semi-definite and positive definite matrices. For a vector $x$,
$X\coloneqq\operatorname{diag}(x)$, $x\circ y$ is the element-wise product, 
$\mathbf{0}_n$ and $\mathbf{1}_n$ are the $n$-dimensional vector of zeros and ones, respectively. Subscript $j$ indexes the problem in the batch and $k$ the interior-point iteration.

\section{Preliminaries: Proximal Interior-Point Method for A Single QP}
\label{sec:prelim}

This section recalls the proximal interior-point method of
\textsc{PIQP}~\cite{piqp} for a single QP. The batch index $j$ is
omitted for convenience.

\subsection{Proximal Method of Multipliers}
\label{subsec:pipm}

We introduce slack variables $s\in\mathbb R^m$ to rewrite the QP as
\begin{equation}
\label{eq:qp_slack}
\begin{aligned}
\min_{x,s}\quad & \tfrac12 x^\top Px + c^\top x \\
\text{s.t.}\quad & Ax=b,\quad Gx-h+s=0,\quad s\ge 0 .
\end{aligned}
\end{equation}
Its augmented Lagrangian with penalty parameter $\delta>0$ is
\begin{equation}
\label{eq:AL}
\begin{aligned}
& \mathcal L^{\mathrm{ALM}}_{\delta}(x,s;\lambda,\nu) \coloneqq
\tfrac12 x^\top Px + c^\top x
+ \lambda^\top(Ax \!-b) \\
& \!+ \tfrac{1}{2\delta}\|Ax \!-\! b\|_2^2
\!+ \nu^\top \!(Gx\! -h \!+s) \!+ \tfrac{1}{2\delta}\|Gx \!-h \!+s\|_2^2,
\end{aligned}
\end{equation}
where $\lambda$ and $\nu$ are the multipliers of the equality and inequality constraints. 
The proximal method of multipliers~\cite{rockafellar1976} iterates
\begin{subequations}
\label{eq:pmm}
\begin{align}
(x^+,s^+) & \in \operatorname*{arg\,min}_{\xi\in\mathbb R^n,\ s\ge 0}\
\mathcal L^{\mathrm{ALM}}_{\delta}(\xi,s;\lambda,\nu) +\! \tfrac{\rho}{2}\|\xi \!-x\|_2^2,
\label{eq:pmm_sub}\\
(\lambda^+,\nu^+) &= (\lambda,\nu) + \tfrac{1}{\delta}\big(Ax^+-b,\ Gx^+-h+s^+\big),
\end{align}
\end{subequations}
with proximal weight $\rho > 0$. One interior-point iteration is applied to solve the
subproblem~\eqref{eq:pmm_sub}. Replacing $s\ge 0$ by a
logarithmic barrier gives, at iteration $k$ with proximal centers
$(\xi_k,\lambda_k,\nu_k)$,
\begin{equation*}
\min_{x,s}\ \mathcal L^{\mathrm{ALM}}_{\delta_k}(x,s;\lambda_k,\nu_k)
- \sigma_k\mu_k\sum_{i=1}^m \ln [s]_i + \tfrac{\rho_k}{2}\|x-\xi_k\|_2^2,
\end{equation*}
where $\mu_k \!> \!0$ is the barrier parameter and $\sigma_k\in(0,1]$ is the
centering parameter defined below. With the auxiliary variables
$y\coloneqq\lambda_k+\frac{1}{\delta_k}(Ax-b)$ and
$z\coloneqq\nu_k+\tfrac1{\delta_k}(Gx-h+s)$, its first-order optimality
conditions read
\begin{subequations}
\label{eq:PMM_IP_KKT}
\begin{align}
Px + c + \rho_k(x-\xi_k) + A^\top y + G^\top z &= 0,\\
Ax - \delta_k(y-\lambda_k) - b &= 0,\\
Gx - \delta_k(z-\nu_k) - h + s &= 0,\\
s\circ z - \sigma_k\mu_k\mathbf 1 &= 0 .
\end{align}
\end{subequations}
Applying Newton's method to~\eqref{eq:PMM_IP_KKT} at the current iterate
$\omega_k=(x_k,y_k,z_k,s_k)$ yields the linear system
\begin{equation} \label{eq:affine}
\underbrace{\begin{bmatrix}
P + \rho_k I & A^\top & G^\top & 0 \\
A & -\delta_k I & 0 & 0 \\
G & 0 & -\delta_k I & I \\
0 & 0 & S_k & Z_k
\end{bmatrix}}_{J_k}
\underbrace{\begin{bmatrix}
\Delta x_k \\ \Delta y_k \\ \Delta z_k \\ \Delta s_k
\end{bmatrix}}_{\Delta\omega_k}
=
\underbrace{\begin{bmatrix}
r^x_k \\ r^y_k \\ r^z_k \\ r^s_k
\end{bmatrix}}_{r_k},
\end{equation}
where $r_k$ is the residual.
The regularization terms $\rho_k I$ and $-\delta_k I$ make $J_k$
nonsingular for any $P\in \mathbb{S}_+^n$ and any $A$, $G$, including the
rank-deficient constraint sets and singular task Jacobians that are common
in practice.


\begin{remark}
When solving the KKT system~\eqref{eq:affine}, we never directly factorize the non-symmetric matrix $J_k$. Eliminating
$\Delta s_k$ yields the symmetric quasi-definite \emph{augmented} system
\begin{equation}
\label{eq:kkt_augmented}
\underbrace{
\begin{bmatrix}
P \!+\!\rho_k I & A^\top & G^\top \\
A & -\delta_k I & 0 \\
G & 0 & -(W_k \!+ \!\delta_k I)
\end{bmatrix}
}_{J_k^\text{aug}}
\begin{bmatrix}\Delta x_k\\ \Delta y_k\\ \Delta z_k\end{bmatrix}
\! = \!
\begin{bmatrix} r_k^x\\ r_k^y\\ \bar r_k^z\end{bmatrix},
\end{equation}
with $W_k=Z_k^{-1}S_k$ and $\bar r_k^z=r_k^z-Z_k^{-1}r_k^s$, which admits an $LDL^\top$ factorization. 
Alternatively, further eliminating $\Delta y_k,\Delta z_k$ gives the \emph{Schur complement} system
\begin{equation}
\label{eq:kkt_complement}
\underbrace{\big( P+\rho_k I+\delta_k^{-1}A^\top A+G^\top(W_k+\delta_k I)^{-1}G \big)}_{J_k^\text{com}}
\Delta x_k=\bar r_k^x,
\end{equation}
where
$\bar r_k^x = r_k^x +\delta_k^{-1}A^\top r_k^y + G^\top(W_k+\delta_k I)^{-1}\bar r_k^z$.
The matrix $J_k^\text{com}\in\mathbb{S}_{++}^n$ admits a Cholesky factorization.
We discuss the choice of factorization strategy in Sec.~\ref{subsec:solver_implementations}.
\end{remark}

\subsection{Mehrotra Predictor--corrector}
\label{subsec:mehrotra}

Each iteration solves~\eqref{eq:affine} twice with the same factorization
of $J_k$~\cite{mehrotra1992}. The \emph{predictor} solves~\eqref{eq:affine}
with $\sigma_k\mu_k=0$ for the affine direction $\Delta\omega^{\mathrm a}_k$
and computes the fraction-to-the-boundary step lengths
\begin{equation}\label{eq:stepsize}
\begin{aligned}
\alpha^{\mathrm a}_p &=\max\{\alpha\in[0,1] \mid s_k+\alpha \Delta s^{\mathrm a}_k\geq (1-\tau)s_k\},\\
\alpha^{\mathrm a}_d &=\max\{\alpha\in[0,1] \mid z_k+\alpha \Delta z^{\mathrm a}_k\geq (1-\tau)z_k\},
\end{aligned}
\end{equation}
with $\tau=0.995$. The centering parameter follows from the predicted
complementarity reduction,
\begin{equation}
\label{eq:parameter}
\begin{aligned}
\eta_k&=\frac{(s_k+\alpha^{\mathrm a}_p\Delta s^{\mathrm a}_k)^\top
     (z_k+\alpha^{\mathrm a}_d\Delta z^{\mathrm a}_k)}{s_k^\top z_k},\\
\sigma_k&=\max\{0,\min\{1,\eta_k\}\}^3 .
\end{aligned}
\end{equation}
The \emph{corrector} sets $\mu_k=s_k^\top z_k/m$ and
solves~\eqref{eq:affine} again with the complementarity residual replaced by
\begin{equation}
\label{eq:rsk}
r^s_k = -s_k\circ z_k - \Delta s^{\mathrm a}_k\circ\Delta z^{\mathrm a}_k + \sigma_k\mu_k\mathbf 1_m,
\end{equation}
giving $\Delta\omega^{\mathrm c}_k$ and step lengths
$(\alpha^{\mathrm c}_p,\alpha^{\mathrm c}_d)$ from~\eqref{eq:stepsize} with
$\Delta\omega^{\mathrm a}_k$ replaced by $\Delta\omega^{\mathrm c}_k$. The
iterate is updated as
\begin{equation}
\label{eq:update}
\begin{aligned}
(x_{k+1},s_{k+1}) &\gets (x_k,s_k)+\alpha^{\mathrm c}_p(\Delta x^{\mathrm c}_k,\Delta s^{\mathrm c}_k),\\
(y_{k+1},z_{k+1}) &\gets (y_k,z_k)+\alpha^{\mathrm c}_d(\Delta y^{\mathrm c}_k,\Delta z^{\mathrm c}_k).
\end{aligned}
\end{equation}

\subsection{Parameter Updates and Termination Criteria}
\label{subsec:peu}

After each step, the proximal centers $(\xi_k,\lambda_k,\nu_k)$ and the
parameters $(\rho_k,\delta_k)$ are updated 
as implemented in
\textsc{PIQP}~\cite[Alg.~2]{piqp}: the dual centers $(\lambda,\nu)$ are moved to
$(y_{k+1},z_{k+1})$ and $\delta$ is reduced aggressively only if the primal
residual has decreased by a fixed factor, and otherwise $\delta$ is reduced
more cautiously; the primal center $\xi$ and $\rho$ are treated in the same
way based on the dual residual. Both parameters are bounded below by
$\underline\delta,\underline\rho>0$ for numerical stability.

The termination is checked according to the absolute and relative residuals as follows:


\begin{equation}
\label{eq:residuals_abs}
\begin{aligned}
r_{\mathrm p} &= \max \{ \left\| Ax-b \right\|_\infty,~ \left\| Gx-h+s \right\|_\infty \},\\
r_{\mathrm d} &= \left\| Px+A^\top y+G^\top z+c \right\|_\infty,\\
r_{\mathrm g} &= \left| x^\top Px+c^\top x+b^\top y+h^\top z \right|,
\end{aligned}
\end{equation}
and their problem-dependent scales
\begin{equation}
\label{eq:residuals_rel}
\begin{aligned}
\gamma_{\mathrm p}&=\max\{\|Ax\|_\infty,\|b\|_\infty,\|Gx\|_\infty,\|h\|_\infty,\|s\|_\infty\},\\
\gamma_{\mathrm d}&=\max\{\|Px\|_\infty,\|A^\top y\|_\infty,\|G^\top z\|_\infty,\|c\|_\infty\},\\
\gamma_{\mathrm g}&=\max\{|x^\top Px|,|c^\top x|,|b^\top y|,|h^\top z|\},
\end{aligned}
\end{equation}
The solver terminates when
\begin{equation}
\label{eq:termination_criteria}
r_{\bullet}\le\max\{\epsilon_{\mathrm{abs}},\ \epsilon_{\mathrm{rel}}\gamma_{\bullet}\},
\qquad \bullet\in\{\mathrm p,\mathrm d,\mathrm g\},
\end{equation}

\section{Differentiation through QPs}
\label{sec:diff}

For learning applications, we differentiate the primal--dual solution
\(\omega(\theta)\coloneqq(x,y,z,s)\) with respect to the QP data
\(\theta\coloneqq(P, A, G, c, b, h)\), aiming at computing $\frac{\partial\omega}{\partial\theta}$.
Let $F(\theta,\omega)$ denote the optimality residual, obtained from~\eqref{eq:PMM_IP_KKT} by setting $\rho=\delta=0$ and the complementarity target $\mu=0$. Together with $s,z\geq0$, the equation $F(\theta,\omega)=0$ gives the KKT conditions for the QP~\eqref{eq:qp_slack}.
At a solution where $\partial F/\partial\omega$ is non-singular, the implicit function theorem gives
\begin{equation}
\frac{\partial\omega}{\partial\theta}=-\left( \frac{\partial F}{\partial\omega} \right)^{-1}\frac{\partial F}{\partial\theta},
\end{equation}
Notice that $\frac{\partial F}{\partial \omega}$ can be approximated by $J_k$ in~\eqref{eq:affine} since they only differ by the regularization terms on the diagonal:
\begin{equation}
    J_k - \frac{\partial F}{\partial \omega} = \operatorname{diag}(\rho_k \mathbf{1}_n, -\delta_k \mathbf{1}_{p},  -\delta_k \mathbf{1}_{m}, \mathbf{0}_{m}),
\end{equation}
where $\rho_k, \delta_k$ are small at the final iteration. Additionally, the regularization terms ensure that $ J_k$ is non-singular. 
Therefore, we approximate the solution sensitivity by
\begin{equation}
\frac{\partial\omega}{\partial\theta}\approx-J_k^{-1} \frac{\partial F}{\partial\theta}.
\end{equation}

Let \(\ell(\omega)\) be a scalar loss that depends on \(\theta\) through \(\omega\). Given the incoming gradient \(\nabla_\omega\ell\), its gradient
w.r.t. $\theta$ is
\begin{equation}
\nabla_\theta\ell
=
\left(\frac{\partial\omega}{\partial\theta}\right)^\top \!\!
\nabla_\omega\ell \\
=
-\left(\frac{\partial F}{\partial\theta}\right)^\top \!
\underbrace{
J_k^{-\top}\nabla_\omega\ell}_{-\eta}
.
\end{equation}
where the adjoint variable \(\eta\) is computed by solving
\begin{equation}
\label{eq:adjoint_system}
J_k^\top\eta=-\nabla_\omega\ell.
\end{equation}

\begin{remark}\label{rmk:reuse_kkt_for_backward}
If factors corresponding to $J_k^\text{aug}$ or $J_k^\text{com}$ at the same $\omega_k, \rho_k, \delta_k$ are available, they can be reused for the adjoint solve~\eqref{eq:adjoint_system}. Otherwise, $J_k$ must be refactorized during the backward pass.
\end{remark}

By reading off the partition \(\eta=(\eta_x,\eta_y,\eta_z,\eta_s)\), the individual data gradients are
\begin{equation}
\label{eq:data_gradients}
\begin{aligned}
\nabla_P\ell
&=\tfrac12(\eta_xx^\top+x\eta_x^\top),
&\quad \nabla_c\ell &=\eta_x,\\
\nabla_A\ell
&=y\eta_x^\top+\eta_yx^\top,
&
\nabla_b\ell
&=-\eta_y,\\
\nabla_G\ell
&=z\eta_x^\top+\eta_zx^\top,
&
\nabla_h\ell
&=-\eta_z.
\end{aligned}
\end{equation}
The expression for \(\nabla_P\ell\) is symmetrized because \(P\) is
parameterized as a symmetric matrix.

\begin{remark}
For sparse matrices with fixed sparsity patterns, the gradients with
respect to \(P\), \(A\), and \(G\) are obtained by gathering the
corresponding dense expressions in~\eqref{eq:data_gradients} only at their
structural nonzeros.
\end{remark}



\section{Batched GPU Execution and Implementation}
\label{sec:method}
\label{sec:impl}

The \cupiqp{} solver executes the proximal interior-point method described in Section~\ref{sec:prelim} through batched operations. 
The solver is implemented in Python and uses NVIDIA \textsc{Warp}~\cite{Macklin_warp} for managing device arrays and generating GPU kernels, with a shared iteration loop and dedicated linear algebra backends for dense, sparse, and multistage problems.

\subsection{GPU-Resident Batch Execution and Termination Control}
\label{sec:data_model}
\label{sec:gpu_execution}

\cupiqp{} is GPU-resident in both numerical data and iteration control.
\footnote{The current sparse backend retains GPU-resident numerical data and computations, but limitations of \textsc{cuDSS} still require a host-side continuation check between iterations.} 
A batch consists of QPs with common dimensions and structural patterns but potentially different numerical data. Problem data and solver variables are stored in device arrays with an explicit batch dimension. All numerical problem data, solver states, and intermediate quantities remain in device memory throughout the solve. The interface accepts device arrays supporting the CUDA array interface through zero-copy views, avoiding intermediate host staging.

Each instance maintains its own primal--dual iterate, proximal centers, regularization parameters, step lengths, and solver status, and advances according to the update rules in Section~\ref{sec:prelim}. Since instances generally require different numbers of iterations, a device-resident activity indicator $a_j\in\{0,1\}$ records whether problem $\mathcal P_j$ requires further iterations. After each iteration, the convergence conditions in~\eqref{eq:termination_criteria} are evaluated independently for each active instance. An instance becomes inactive upon convergence or detection of another terminal solver status, including unboundedness and infeasibility. Subsequent updates to its iterate, proximal centers, and solver parameters are disabled, preserving its final solution and termination status. 

Step lengths, parameter updates, termination checks, and activity updates are computed on the GPU. The activity indicators are reduced on the device to form the aggregate continuation flag $a^{\mathrm{batch}}=\bigvee_{j=1}^{B}a_j$. The shared loop terminates when this flag is false, or the prescribed iteration limit is reached. Both per-instance numerical decisions and batch-level loop control therefore remain on the device, avoiding per-iteration host intervention and host--device synchronization unless verbose output is enabled to print solver progress at each iteration.

\subsection{Execution Overhead Reduction}
\label{subsec:performance_optimization}

For batches of small QPs, kernel dispatch and memory allocation can account for a substantial fraction of the solve time. \cupiqp{} makes efforts to reduce these costs through kernel fusion, workspace reuse, and CUDA graph replay.

Operations for residual evaluation, variable updates, step-length computation, regularization updates, and termination checks are grouped into a small number of \textsc{Warp} kernels, reducing kernel launches and intermediate memory traffic. The required reductions are performed independently within each QP using per-problem block reductions.

During the \texttt{setup} phase, \cupiqp{} allocates reusable solver workspaces, generates and compiles the required \textsc{Warp} kernels, and performs structure-dependent preprocessing like re-ordering and symbolic analysis where required. The allocated workspaces, compiled kernels, and preprocessing results are reused across subsequent \texttt{update}, \texttt{solve}, and optional \texttt{backward} phases, keeping allocation, compilation, and preprocessing outside repeated solves. The persistent workspaces also provide stable memory addresses for CUDA graph capture and replay. 

The device-side iteration control described in Section~\ref{sec:gpu_execution} is implemented using CUDA conditional graph nodes to determine whether to execute another iteration.
A complete iteration is captured once and replayed for subsequent iterations, reducing the host-side overhead of dispatching individual kernels and library calls. Verbose execution instead uses a host-controlled loop to report progress at each iteration.
The interface also supports capturing \texttt{update} and \texttt{solve} within a user-defined CUDA graph when full-solve capture is supported and verbose output is disabled. This allows simulation, QP construction, solution, and actuation to form a single replayable pipeline in robotic applications.

\subsection{Problem Representations and Linear Algebra Backends}
\label{subsec:solver_implementations}

\cupiqp{} provides dedicated interfaces for dense, sparse, and multistage
QPs. All three follow the batched algorithm of Section~\ref{sec:method} and
share the surrounding implementation, but differ in input representation,
data layout, and KKT solution strategy.

\subsubsection{Dense solver}
All problem data are given as dense arrays and stored, together with all solver variables, as contiguous device arrays with a leading batch dimension, so that every operation of the iteration is a uniform batched operation without indirect memory access. The dense solver uses the Schur
complement system~\eqref{eq:kkt_complement}, which is factorized and solved with the batched Cholesky routines. The matrix $A_j^\top A_j$ is computed once and cached for reuse across iterations.
For small matrices (empirically smaller than $32 \times 32$), we use shape-specialized kernels for linear algebra operations, while using \textsc{cuBLAS} and \textsc{cuSOLVER} for larger matrices.

\subsubsection{Sparse solver}
The matrices $P_j$, $A_j$, and $G_j$ are given in compressed sparse rows (CSR) format with one sparsity pattern for the whole batch. The index arrays are stored once, and the nonzero values are stored as contiguous batched device arrays. The sparse solver uses the augmented system~\eqref{eq:kkt_augmented} and performs $LDL^\top$ factorization through \textsc{cuDSS}. The symbolic analysis is performed once in the \texttt{setup} phase and reused for the subsequent numerical factorization and triangular solve.

\subsubsection{Multistage solver}
The multistage solver targets batched QPs arising from multistage optimization problems following~\cite{piqp_multistage}, such as optimal control problems (OCPs). For $j=1,\ldots,B$, each $\mathcal{P}_j$ has the form
\begin{equation*}
\begin{aligned}
\min_{x_j} \quad
& \sum_{i=0}^{N-1}
\left(
\frac{1}{2}
\begin{bmatrix}
x_{j,i} \\
x_{j,i+1}
\end{bmatrix}^{\!\top}
\begin{bmatrix}
\bar{Q}_{j,i} & \bar{S}_{j,i}^{\top} \\
\bar{S}_{j,i} & 0
\end{bmatrix}
\begin{bmatrix}
x_{j,i} \\
x_{j,i+1}
\end{bmatrix}
+ \bar{c}_{j,i}^{\top} x_{j,i}
\right)
\\
& \qquad
+ \frac{1}{2}x_{j,N}^{\top}\bar{Q}_{j,N}x_{j,N}
+ \bar{c}_{j,N}^{\top}x_{j,N}
\\
\text{s.t.}\quad
& \bar{A}_{j,i}x_{j,i}
+ \bar{B}_{j,i}x_{j,i+1}
= \bar{b}_{j,i},
\;\qquad i=0,\ldots,N\!-\!1,
\\
& \bar{C}_{j,i}x_{j,i}
+ \bar{D}_{j,i}x_{j,i+1}
\leq \bar{h}_{j,i},
\qquad i=0,\ldots,N \!-\!1,
\\
& \bar{A}_{j,N}x_{j,N}=\bar{b}_{j,N}, \\
& \bar{C}_{j,N}x_{j,N}\leq \bar{h}_{j,N},
\end{aligned}
\end{equation*}
where $x_{j,i}\in\mathbb{R}^{n_i}$ denotes the decision variable at stage $i$ of $\mathcal{P}_j$, and $N$ is the horizon length. 
We assume uniform stage sizes, with $p_i$ equality constraints and $m_i$ inequality constraints at every stage.\footnote{For OCPs, the terminal stage typically has no control input and therefore fewer variables. We can pad it to obtain a uniform stage size.}
All QPs in the batch share $N$, $n_i$, $p_i$, $m_i$, while their numerical data may vary with $j$.
The stage blocks are given directly and stored as batched dense arrays and never assembled into big sparse matrices. 

With the Schur complement formulation~\eqref{eq:kkt_complement}, the matrix $J_k^\text{com}$ has a block-tri-diagonal structure. We solve it using the parallel method proposed in~\cite{schwan2026socu} and its open-source implementation \textsc{socu}, which achieves an $\mathcal O(\log N)$ time complexity
\footnote{The $\mathcal{O}(\log N)$ time complexity assumes sufficient hardware resources to execute
the operations concurrently. With a sufficiently large batch size, the GPU can saturate the available parallelism and the wall-clock time no longer exhibits logarithmic scaling w.r.t. $N$.}
thanks to the parallelism exploited both across the batch and along the horizon, while the sequential approaches in~\cite{piqp_multistage} have 
$\mathcal O(N)$ time complexity.

\subsection{Constraint Handling and Numerical Robustness}
\label{subsec:numerical_implementation}

Two-sided inequalities $h_l\leq Gx\leq h_u$ and variable bounds
$x_l\leq x\leq x_u$ arise frequently in robotics. 
\cupiqp{} accepts both forms directly and incorporates their contributions using the treatment
in~\cite[Ch.~3]{schwan_thesis}, avoiding enlarging the Schur complement system~\eqref{eq:kkt_complement}.

To improve numerical conditioning, \cupiqp{} applies Ruiz equilibration~\cite{ruiz2001scaling}. Optional iterative refinement is available after each linear solve. The maximum number of refinement iterations is user-configurable. 

\subsection{Implicit Differentiation}
\label{subsec:impl_diff}

We implement the implicit differentiation method in Section~\ref{sec:diff}.
When factors of the KKT system at the differentiation point are available, they can be reused for the adjoint solve without an additional factorization, as mentioned in Remark~\ref{rmk:reuse_kkt_for_backward}.
However, when differentiating through closed-loop trajectories generated by an optimization-based policy, as demonstrated in Sec.~\ref{subsec:exp_diff_mpc}, retaining the factors from every forward solve can require substantial memory.
We therefore only save the problem coefficients and solutions needed to reconstruct each system and perform numerical factorization again during the backward pass, trading additional computation for reduced memory usage.
We also implement dedicated \textsc{PyTorch} and \textsc{JAX} interfaces to integrate the forward solve and backward computation into their respective automatic differentiation pipelines.





\begin{figure*}[t]
\centering
\includegraphics[width=1.0\textwidth]{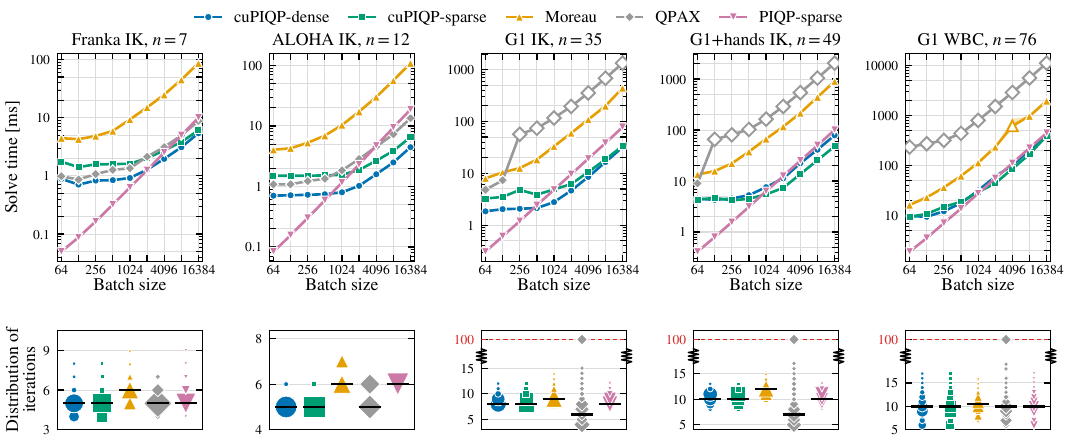}
\caption{Batched solve time (top) and distribution of iterations at $B=16384$ (bottom) for the IK and  WBC problems, with the median marked by solid black bars. The iteration limit is marked with a red dashed line. The hollow markers indicate that for at least one QP in the batch, the solver failed to converge to the required tolerance within the maximum number of iterations.}
\label{fig:ik_scaling}
\end{figure*}

\section{Experimental results}
\label{sec:exp}

We comprehensively evaluate the \textsc{cuPIQP} solver on batched QPs from multiple robotic applications,
including inverse kinematics, whole-body control, and differentiable linear MPC.

\subsection{Experimental Setting}
All experiments run on Ubuntu 24.04 with an Intel Ultra 9 285K CPU, 64 GB of RAM, and an NVIDIA GeForce RTX 5090 (32 GB), under \texttt{float64} arithmetic. 
We compare \textsc{cuPIQP} against \textsc{QPAX} and \textsc{Moreau}, with the latter using \textsc{cuDSS} as its linear-system solver.  Iterative refinement is disabled if applicable.
We additionally include \textsc{PIQP}~\cite{piqp}'s sparse backend and \textsc{HPIPM}\footnote{We use its \texttt{balance} mode, same as the default setting in \textsc{acados}.}~\cite{frison2020hpipm} as CPU solver baselines, using custom C++ bindings to create 8 solver instances pinned onto 8 performance cores of the CPU to solve the batched QPs concurrently with OpenMP. 

For all solvers, we use absolute tolerance $10^{-6}$ and relative tolerance $10^{-7}$, except \textsc{QPAX} because it struggles to converge to high accuracy and frequently runs to the maximum iteration. Therefore, we relax its tolerance to $10^{-4}$, giving it a potential
advantage in both iteration count and solve time.

We report steady-state wall time with solver data already resident in host or device memory. Data preparation, setup, JIT compilation, and CUDA-graph capture are excluded. The first call warms up the solvers under every configuration. 
Reported times are medians of at least 20 solves in Sec.~\ref{subsec:exp_ik_wbc} and 10 solves in Sec.~\ref{subsec:exp_diff_mpc}.

\subsection{Batched Inverse Kinematics and Whole-body Control}
\label{subsec:exp_ik_wbc}

\subsubsection{Inverse kinematics}
We form IK QPs such that each QP is a linearized SQP step of the nonlinear IK problem at the current configuration $q$.
Given the Jacobian $J_i(q)$ and pose error $d_i(q)$ of each end-effector $i$, we compute a joint increment $\Delta q$ by weighting its task-space tracking residual with $W_i\succeq0$. We set the damping weight to $w=10^{-3}$ and use $\gamma=0.95$ to retain a margin from the joint limits:
\begin{equation*}
\begin{aligned}
\min_{\Delta q} \quad & \tfrac12 \sum\nolimits_i\|J_i(q)\,\Delta q - d_i(q)\|_{W_i}^2 + \tfrac{1}{2}w \|\Delta q\|^2\\
\text{s.t.} \quad & \gamma\,(q_{\min}-q)\le\Delta q\le\gamma\,(q_{\max}-q).
\end{aligned}
\end{equation*}

The two Unitree G1 benchmark datasets record the QPs generated by GMR~\cite{araujo2025gmr} while retargeting all 77 motion-capture sequences from the LAFAN1 dataset~\cite{harvey2020lafan1} that span diverse motions including walking, running, and fighting. 
Each recorded QP is sampled from a different motion frame, while the dataset collectively covers different SQP iterations.
Together, these variations yield a diverse benchmark dataset that captures both the changing robot configurations and the different numerical regimes encountered throughout SQP convergence in practical applications.
The benchmark datasets for Franka and ALOHA robots are constructed in the same manner using real-world teleoperation data~\cite{khazatsky2024droid} and ~\cite{zhao2023aloha}, respectively.



\subsubsection{Batched whole-body control}

We form inverse dynamics WBC QPs for the 29-DoF Unitree G1, generated from LAFAN1 motions retargeted to the robot~\cite{harvey2020lafan1}. Following a standard inverse dynamics formulation~\cite{kuindersma2016optimization}, each QP has 76 variables comprising generalized accelerations, joint torques, and two 6D sole wrenches; 47 rigid-body-dynamics and contact-acceleration equalities; 20 friction, center-of-pressure, and torsional inequalities; and joint,
torque, acceleration, and wrench bounds. Left support, right support, flight,
and double support are represented with a common dimension and sparsity
pattern by setting airborne contact wrenches to zero. The objective tracks
the recorded whole-body motion and regularizes torques and contact wrenches.

\subsubsection{Results}
We compare \cupiqp{}-dense, \cupiqp{}-sparse, \textsc{Moreau}, \textsc{QPAX}, and \textsc{PIQP} for $B=64, 128, \ldots,16384$. Recomputed constraint violations are below $10^{-6}$ and objectives agree with each other within $10^{-5}$, except for \textsc{QPAX}, which reaches the iteration limit without converging on some problems. At $B=16384$, \textsc{QPAX} fails on 56 G1 problems and 215 G1-with-hands problems.

Fig.~\ref{fig:ik_scaling} reports the batch-scaling solver time and the distribution of iterations for problems in the largest batch. 
Among the dense solvers, \cupiqp{}-dense consistently achieves lower solve times than \textsc{QPAX}, despite using a stricter stopping tolerance. The difference is most pronounced on the G1 IK and WBC problems, where \textsc{QPAX} reaches the iteration limit for a subset of instances.
Among the sparse solvers, \cupiqp{}-sparse consistently outperforms \textsc{Moreau} with up to $18.8\times$ speedup. 

Importantly, the speedups are not obtained by reducing the number of interior-point iterations: the iteration distributions of both \cupiqp{} backends exhibit tightly clustered iteration counts, with medians comparable to those of \textsc{Moreau} and \textsc{PIQP}. Therefore, the substantial timing advantage over \textsc{Moreau} and \textsc{QPAX} is consistent with the emphasis on performance optimization in Section~\ref{subsec:performance_optimization}.

The \textsc{cuPIQP} dense solver remains preferable for the small problems, whereas its sparse solver becomes about $1.65\times$ faster than the dense solver for the G1-with-hands problems, demonstrating the benefit of exploiting sparsity once the problem is sufficiently large. 

Compared to the CPU solver \textsc{PIQP}, 8-core \textsc{PIQP} is faster at small batch sizes across all five benchmarks, but is surpassed by \cupiqp{} as the batch size increases due to the improved GPU utilization and amortization of fixed execution costs. 
At the largest batch, \cupiqp{} is faster than \textsc{PIQP} by $1.81\times$ on Franka IK, $2.47\times$ on G1 IK, $1.29\times$ on G1--with--hands IK, and $4.26\times$ on ALOHA. 


\subsection{Learning Differentiable Linear MPC for Quadrotors}
\label{subsec:exp_diff_mpc}

\begin{figure*}[t]
\centering
\includegraphics[width=0.95\textwidth]{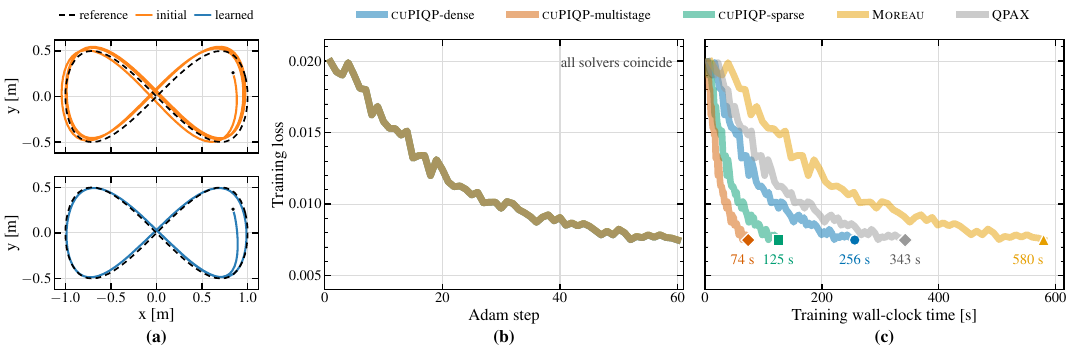}
\caption{Learning the MPC weights at $B=256$. (a)~Top-down paths of a closed-loop rollout with the initial weights (top) and with the learned weights (bottom). The dot marks the starting position. (b)~Training loss versus Adam step. (c)~Loss versus wall-clock training time; markers show the end of the 60 steps with the training time.}
\label{fig:learn_curves}
\end{figure*}

\begin{table*}[t]
\centering
\caption{Time per Adam step (seconds) at multiple batch sizes, and its percentage share of each component, including QP solves for the forward rollout, QP backward passes, and everything else. For HPIPM, only its forward QP solve time is reported. 
Shaded: fastest QP forward time per batch size among the GPU solvers.}
\label{tab:learn_breakdown}
\footnotesize
\setlength{\tabcolsep}{2.1pt}
\setlength{\fboxsep}{1pt}
\renewcommand{\arraystretch}{1.12}
\begin{tabular}{@{}l|cccc|cccc|cccc@{}}
\toprule
& \multicolumn{4}{c|}{$B=64$} & \multicolumn{4}{c|}{$B=256$} & \multicolumn{4}{c}{$B=1024$} \\
\cmidrule(lr){2-5}\cmidrule(lr){6-9}\cmidrule(lr){10-13}
Solver & QP fwd. & QP bwd. & Other & Total & QP fwd. & QP bwd. & Other & Total & QP fwd. & QP bwd. & Other & Total \\
\hline
\cupiqp{}-dense & 1.37 (80\%) & 0.09 (5\%) & 0.26 (15\%) & 1.72 & 3.83 (88\%) & 0.25 (6\%) & 0.27 (6\%) & 4.36 & 14.26 (92\%) & 0.92 (6\%) & 0.29 (2\%) & 15.47 \\
\cupiqp{}-multistage & \colorbox{cyan!30}{0.43} (60\%) & 0.04 (5\%) & 0.25 (35\%) & 0.72 & \colorbox{cyan!30}{0.91} (74\%) & 0.07 (6\%) & 0.25 (20\%) & 1.24 & \colorbox{cyan!30}{2.88} (87\%) & 0.18 (5\%) & 0.26 (8\%) & 3.31 \\
\cupiqp{}-sparse & 0.72 (70\%) & 0.06 (5\%) & 0.25 (24\%) & 1.03 & 1.74 (82\%) & 0.12 (6\%) & 0.26 (12\%) & 2.12 & 5.28 (90\%) & 0.30 (5\%) & 0.27 (5\%) & 5.85 \\
\textsc{Moreau} & 2.89 (82\%) & 0.38 (11\%) & 0.26 (7\%) & 3.53 & 8.79 (86\%) & 1.16 (11\%) & 0.27 (3\%) & 10.21 & 35.18 (89\%) & 4.20 (11\%) & 0.28 (1\%) & 39.66 \\
\textsc{QPAX} & 2.27 (85\%) & 0.12 (5\%) & 0.27 (10\%) & 2.66 & 5.14 (90\%) & 0.24 (4\%) & 0.34 (6\%) & 5.72 & 18.03 (94\%) & 0.74 (4\%) & 0.34 (2\%) & 19.11 \\
\hline
\textsc{HPIPM} & 0.21 & -- & -- & -- & 0.79 & -- & -- & -- & 3.13 & -- & -- & -- \\
\bottomrule
\end{tabular}
\end{table*}



To evaluate \cupiqp{} in the end-to-end differentiable learning pipeline, we optimize the cost weights of a linear MPC for a Crazyflie quadrotor tracking a figure-eight trajectory with a period of \SI{4}{s}.

\subsubsection{MPC formulation}
The MPC uses a hover-linearized model with state
\(x=(p,\phi,v,\omega)\in\mathbb{R}^{12}\), comprising position, roll--pitch--yaw angles, linear velocity, and body-frame angular velocity. The input \(u\in\mathbb{R}^{4}\) contains the four rotor thrusts.
The prediction horizon of \SI{1}{s} is discretized into \(N=20\) steps, with dynamics
\(x_{k+1}=\hat A x_k+\hat B u_k + \hat b\) and initial state \(x_0=\hat x\).
With \(e_k:=x_k-x_k^{\mathrm{ref}}\), the objective is
\[
\tfrac12\sum\nolimits_{k=0}^{N-1}
\left(\|e_k\|_Q^2+\|u_k\|_R^2\right)
+\tfrac12\|e_N\|_Q^2.
\]
Roll and pitch are bounded by \(\pm\SI{0.45}{rad}\), 
body-frame angular velocities by \(\pm \SI{5}{rad/s}\), 
and each rotor thrust by \([\SI{0}{N}, \SI{0.149}{N}]\).
To ensure feasibility, all state constraints are softened with linear and quadratic slack penalties with
coefficients fixed throughout training. The input constraints remain hard.

We fix \(R \!=\!I_4\) and parameterize the diagonal state weights as
$Q(\varphi)=\operatorname{diag}(q_1,\ldots,q_{12})$ with 
$q_i=10^{-3}+\epsilon_i\operatorname{softplus}(\varphi_i)$, where $\epsilon_i>0$ are fixed scales. 
We deliberately initialize $\varphi$ to yield poor tracking performance.
\textsc{QPAX} and \cupiqp{}-dense solve the condensed QP obtained by eliminating the state variables, and the other solvers retain both states and inputs as decision variables.

\subsubsection{Learning problem}
The MPC weights are trained to minimize a closed-loop objective defined by the per-step loss
\begin{equation*}
\begin{aligned}
\ell_k= &
\|p_k \!-\! p_k^{\mathrm{ref}}\|_2^2
\! + \! 0.1\|v_k \!-\! v_k^{\mathrm{ref}}\|_2^2
\! + \! 0.1\|\phi_k^{\mathrm{rp}} \!-\! \phi_k^{\mathrm{ref}}\|_2^2 \\
& +0.01\|u_k\|_2^2
+0.05\|u_k \!-\! u_{k-1}\|_2^2
+100\,\mathrm{viol}_k ,
\end{aligned}
\end{equation*}
where $\phi_k^{\mathrm{rp}}$ contains roll and pitch,
and $\mathrm{viol}_k$ is the sum of squared violations of the tilt, angular-rate, and altitude limits.
The training objective averages this loss over $T=40$ rollout steps and $B=256$ trajectories:
\[
\mathcal L(\varphi) = \tfrac{1}{BT} \sum\nolimits_{j=1}^{B}\sum\nolimits_{k=0}^{T-1}\ell_{j,k}.
\]

We optimize $\varphi$ with 60 Adam steps. Each step evaluates a batch of \SI{2}{s} closed-loop rollouts using a nonlinear quadrotor simulation model, with randomized initial conditions and figure-eight dimensions. 
The controller solves a new MPC problem at each rollout step, yielding 40 batched QP solves per Adam step. 
Gradients propagate through these solves and the nonlinear model in the simulator to the entries of $Q$.
All solvers use identical training batches, initial weights, loss functions, and Adam optimizer settings.

\subsubsection{Results}
We evaluate the controllers on 64 test trajectories not used during training.
Each closed-loop simulation lasts \SI{12}{s} and covers three complete cycles of the figure-eight reference.
During training, all QPs are solved to the specified tolerance.
After training, the root-mean-square (RMS) position error decreases from \SI{187}{mm} to \SI{55}{mm} with the learned weights for all five solvers, as illustrated by Fig.~\ref{fig:learn_curves}(a).
Moreover, the loss curves versus Adam step nearly coincide across solvers as shown in Fig.~\ref{fig:learn_curves}(b).
Despite differences in QP formulation and differentiation implementation, the solvers therefore produce essentially the same learning progress. This agreement is consistent with numerically consistent forward solutions and backward gradients in this experiment.
The solvers therefore achieve essentially the same control performance and differ primarily in computational time, as shown in Fig.~\ref{fig:learn_curves}(c).

On the condensed formulation, \cupiqp{}-dense reduces training time by $25\%$ relative to \textsc{QPAX}, even under stricter convergence requirements.
On the uncondensed formulation, \cupiqp{}-sparse provides a $4.6\times$ speedup over \textsc{Moreau}.
The \cupiqp{}-multistage provides a further $1.7\times$ speedup over \cupiqp{}-sparse and completes training $4.6\times$ and $7.8\times$ faster than \textsc{QPAX} and \textsc{Moreau}, respectively, showing the advantages of explicitly exploiting the multistage structure in optimal control problems.

While Fig.~\ref{fig:learn_curves} reports training at $B=256$, Table~\ref{tab:learn_breakdown} additionally
evaluates $B=64$ and $B=1024$ to examine how the per-step computational cost scales with batch size and breaks down the total time into QP forward solves, QP backward passes, and the remaining computation.
Forward QP solves account for $60$--$94\%$ of the total time across the differentiable solvers and batch sizes, whereas backward passes account for only approximately $5$--$6\%$ for \cupiqp{}.
The multistage backend achieves the lowest total time at all three batch sizes, with its relative
advantage increasing as the batch grows. At $B=1024$, an Adam step takes \SI{3.31}{s}, giving speedups of $1.8\times$ and $4.7\times$ over the sparse and dense \cupiqp{} backends, and $5.8\times$ and $12.0\times$ over \textsc{QPAX} and \textsc{Moreau}, respectively.

We additionally report the forward QP solve time of \textsc{HPIPM} in Table~\ref{tab:learn_breakdown} since it does not support differentiation. \textsc{HPIPM} is faster at $B=64$ and $B=256$, while \cupiqp{}-multistage
becomes competitive at $B=1024$.  At this batch size, its forward time is \SI{2.88}{s}, compared with \SI{3.13}{s} for \textsc{HPIPM}.




\section{Conclusions}
\label{sec:conclusions}

We presented \cupiqp{}, a batched and differentiable QP solver on GPUs that supports dense, sparse, and multistage problems within a unified framework. Benchmark results demonstrate solve-time speedups of up to $18.8\times$ over \textsc{Moreau} and $4.26\times$ over eight-core \textsc{PIQP}. In differentiable quadrotor MPC, these computational gains translated into end-to-end training speedups of $4.6\times$ over \textsc{QPAX} and $7.8\times$ over \textsc{Moreau}.
By reducing the computational cost of integrating model-based optimization into large-scale robot learning and control, \cupiqp{} provides a reusable foundation for future research and applications.

Future work will extend this foundation to differentiable sequential quadratic programming for nonlinear MPC and broader optimization-based applications.

\section*{Acknowledgment}
Anthropic Claude Code and OpenAI Codex were used to assist with manuscript preparation, software debugging and performance optimization, and the implementation of experimental scripts.
All generated material was reviewed, verified, and revised by the authors, who take full responsibility for the content and results.

\balance
\bibliographystyle{IEEEtran}
\bibliography{reference}

\end{document}